\documentclass[12pt,leqno]{amsart}

\usepackage[utf8]{inputenc}   
\usepackage[T1]{fontenc}      
\usepackage{geometry}         
\usepackage[francais]{babel}  

\usepackage{amssymb}
\usepackage{enumerate}

\voffset-1cm \markleft{\rm Pfaffian Systems} 
\newtheorem{theorem}{Theorem}
\newtheorem{definition}[theorem]{Definition}
\newtheorem{corollary}[theorem]{Corollary}

\newtheorem{lemma}[theorem]{Lemma}

\newtheorem{proposition}[theorem]{Proposition}

\newtheorem{defi}{D\'{e}finition}

\newcommand{\K}{\mathbb K}

\newcommand{\g}{\frak{g}}

\newcommand{\R}{\mathbb R}

\newcommand{\ra}{\rightarrow}

\newcommand{\pf}{\noindent{\it Proof. }}
\newcommand{\ds}{\displaystyle}

\newcommand{\we}{\wedge }

\newcommand{\h}{\mathfrak h}

\usepackage{graphicx}

\title{Contact forms on the simple Lie Group $SL(2p,\K)$}
\author{Elisabeth Remm}
\date{}
\address{Universit\'e de Haute-Alsace, IRIMAS UR 7499, F-68100 Mulhouse, France.}
\email{elisabeth.remm@uha.fr}

\begin{document}

\maketitle

\section{Introduction: Contact manifolds}
A contact form on a $(2p+1)$-dimensional differential manifold is a Pfaffian form $\alpha$ such that
$\alpha \we (d\alpha)^p(x) \neq 0$ for all $x \in M$.  A contact structure is a $2p$-dimensional distribution defined by a contact form.  Two contact structures are isomorphic when the contact forms which define these structures are conjugated by a diffeommorphism $h$, that is $h^*\alpha \we \alpha' =0.$ The construction of contact manifolds, or the existence of contact forms on a manifold is a major problem in differential geometry. For example, in mechanics, the importance of dynamic systems admitting for integral invariant a form of contact was demonstrated by G.Reeb. Here are some major results in this direction: Any $3$-orientable $3$-closed manifold is a contact manifold. In dimensions greater than $3$, the results
are few: we know contact structures on the
nodes associated with complex hypersurfaces, especially
exotic spheres bordering a parallel variety and
some principal fiber bundles.  The first examples of "non classical"contact structure are described by Boothby and Wang.  These structures of Boothby-Wang 
are both invariant and transverse to the orbits of a free action
of $S^1$. Starting of this point of view, Lutz introduces the concept of contact structures invariant by an action of a group on the manifold. He shows that exists contact structures on the $5$-dimensional torus (for the $1$ and $3$-dimensional cases, it is not difficult to presents such structures). Latter, Bourgeois generalizes the Lutz result for all  odd dimensional torus. In the same time, in the same context, Goze studied the left invariant contact form on a Lie group. For example, he proves that there exists such contact forme on  a semi-simple Lie group if and only its the rank is equal to $1$ and he introduces the notion of degree of freedom  of a contact structure. To end this presentation, which is very not exhaustive, the study of contact form on  the connected sum of two contact manifolds by Meckert \cite{Mecker} or the study of the set of contact forms on $3$-dimensional manifolds by Hadjar \cite{Hadjar}.

The aim of the present work is to look if  simple Lie groups are contact manifold. Since left invariant form cannot be contact as soon as the rank of the Lie group is better than $2$, we will study non invariant form on such Lie group. We begin this work by the study of $SL(n)$. 

To conclude this introduction, let us recall that the problem of the existence of invariant contact forms on non-simple Lie groups has been the subject of numerous works, let us cite for example \cite{Sa1} where contact forms are constructed on nilpotent Lie groups by a extension process and \cite{Sa2} where one associates invariant such a Pfaffian to any left invariant contact forms on a Lie group.

\section{On the Cartan class of left invariant pfaffian forms on a simple Lie group}
Let $M$ be a $n$-dimensional differentiable manifold and $\alpha$ a Pfaffian form on $M$, that is a differential form of degree $1$. The
{\it characteristic space of $\alpha$ at a point $x \in M$} is the linear subspace~$\mathcal{C}_x(\alpha)$ of the tangent space $T_x(M)$ of $M$ at the point $x$ defined by
\begin{gather*} \mathcal{C}_x(\alpha)=A(\alpha(x)) \cap A({\rm d}\alpha(x)),\end{gather*}
where \begin{gather*} A(\alpha(x))=\{X_x \in T_x(M), \, \alpha(x)(X_x)=0\}\end{gather*} is the {\it associated subspace of $\alpha (x)$},
\begin{gather*} A({\rm d}\alpha(x))= \{X_x \in T_x(M), \, {\rm i} (X_x){\rm d} \alpha(x) =0 \}\end{gather*} is the {\it associated subspace of ${\rm d}\alpha (x)$} and ${\rm i}(X_x){\rm d}\alpha(x)(Y_x)={\rm d}\alpha(x)(X_x,Y_x)$.

\begin{definition}Let $\alpha$ be a Pfaffian form on the differential manifold $M$. The {\it Cartan class of~$\alpha$ at the point $x \in M$}~\cite{Cartan} is the codimension of the characteristic space~$\mathcal{C}_x(\alpha)$ in the tangent space~$T_x(M)$ to~$M$ at the point~$x$. We denote it by $\operatorname{cl}(\alpha)(x)$.
\end{definition}

 The function $x \rightarrow \operatorname{cl}(\alpha)(x)$ is with positive integer values and is lower semi-continuous, that is, for every $x \in M$ there exists a suitable neighborhood $V$ such that for every $x_1 \in V$ one has $\operatorname{cl}(\alpha)(x_1) \geq \operatorname{cl}(\alpha)(x)$.

The {\it characteristic system of $\alpha$ at the point $x$ of $ M$} is the subspace $ \mathcal{C}_x^*(\alpha)$ of the dual $T_x^*(M)$ of $T_x(M)$ orthogonal to $\mathcal{C}_x(\alpha)$:
\begin{gather*} \mathcal{C}_x^*(\alpha)=\{\omega(x) \in T_x^*(M), \, \omega(x)(X_x)=0, \, \forall\, X_x \in \mathcal{C}_x(\alpha)\}.\end{gather*}
Then
\begin{gather*}\operatorname{cl}(\alpha)(x)= \dim \mathcal{C}_x^*(\alpha).\end{gather*}

\begin{proposition}
If $\alpha$ is a Pfaffian form on $M$, then
\begin{itemize}\itemsep=0pt
\item $\operatorname{cl}(\alpha)(x)= 2p+1$ if $(\alpha \wedge ({\rm d}\alpha)^p)(x) \neq 0$ and $ ({\rm d}\alpha)^{p+1}(x)=0$,
\item $\operatorname{cl}(\alpha)(x)=2p$ if $({\rm d}\alpha )^p(x) \neq 0 $ and $ (\alpha \wedge ({\rm d}\alpha )^{p})(x)=0$.
\end{itemize}
\end{proposition}
In the first case, there exists a basis $\{\omega_1(x)=\alpha(x),\omega_2(x),\dots,\omega_n(x)\}$ of $T_x^*(M)$ such that
\begin{gather*} {\rm d}\alpha(x)=\omega_2(x)\wedge \omega_3(x)+ \dots +\omega_{2p}(x)\wedge \omega_{2p+1}(x)\end{gather*}
and
\begin{gather*} \mathcal{C}_x^*(\alpha)=\R\{\alpha(x)\}+A^*({\rm d}\alpha (x)).\end{gather*}
In the second case, there exists a basis $\{\omega_1(x)=\alpha(x),\omega_2(x),\dots,\omega_n(x)\}$ of $T_x^*(M)$ such that
\begin{gather*} {\rm d}\alpha(x)=\alpha(x)\wedge \omega_2(x)+ \dots +\omega_{2p-1}(x)\wedge \omega_{2p}(x)\end{gather*}
and
\begin{gather*}\mathcal{C}_x^*(\alpha)=A^*({\rm d}\alpha (x)).\end{gather*}

If the function $\operatorname{cl}(\alpha)(x)$ is constant, that is, $\operatorname{cl}(\alpha)(x)= \operatorname{cl}(\alpha)(y)$ for any $x, y \in M$, we say that the Pfaffian form $\alpha$ is of constant class and we denote by $\operatorname{cl}(\alpha)$ this constant. The distribution
\begin{gather*} x \rightarrow \mathcal{C}_x(\alpha)\end{gather*}
is then regular and it is an integrable distribution of dimension $n-\operatorname{cl}(\alpha)$, called the characteristic distribution of~$\alpha$. It is equivalent to say that the Pfaffian system
\begin{gather*} x \rightarrow \mathcal{C}^*_x(\alpha)\end{gather*}
is integrable and of dimension $\operatorname{cl}(\alpha)$.

Assume now that $M=G$ is a Lie group. We identify any left invariant Pfaffian form $\alpha$ on a Lie group $G$ with the corresponding element, denoted also $\alpha$, of the dual $\g^*$ of the Lie algebra $\g$ of $G$. Modulo this identification, we have
$$d\alpha (X,Y)=-\alpha[X,Y]$$
for any $X,Y \in \g$.
The {\it Cartan class $\operatorname{cl}(\alpha)$ of the linear form $\alpha \in \g^*$} is defined as follow:
\begin{itemize}\itemsep=0pt
\item $\operatorname{cl}(\alpha)=2p+1$ if and only if $\alpha \wedge ({\rm d}\alpha)^p\neq 0$ and $ ({\rm d}\alpha)^{p+1}=0$,
\item $\operatorname{cl}(\alpha)=2p$ if and only if $({\rm d}\alpha )^p \neq 0 $ and $ \alpha \wedge ({\rm d}\alpha )^{p}=0$.
\end{itemize}
For example, if $\g$ is a nilpotent Lie algebra or if it is a compact real Lie algebra, then class of any non trivial linear form is odd. 

If there is on a $(2p+1)$-dimensional Le algebra a linear form of class $2p+1$, that is a linear contact form, then $\g$ is called a contact Lie algebra. In this case $G$ is a contact manifold. The Heisenberg Lie group is a contact manifold. One can characterize the set of $(2p+1)$ dimensional contact Lie algebra in term of deformations of the Heisenberg algebra (\cite{G.R.JDGA}). In the simple or semi-simple case we have:
\begin{proposition}\cite{BW,GozeCras1}\label{simple}
Let $\g$ be a real or complex semi-simple Lie algebra of rank $r$. Then, for any $\omega \in \g^*$ we have
$$cl(\omega) \leq n-r+1.$$
\end{proposition}
\pf For any $\omega  \in g^*$ with $\omega  \neq 0$, the associated space $A(d\omega)=\{X \in \g, i(X)d\omega=0\}$ is a subalgebra of dimension less than or equal to $n-c+1$ where $n=\dim \g$ and $c$ is the Cartan class of $\omega$. If $K$ denotes the Killing-Cartan form on $\g$, since $K$ is non degenerate, there exists $U_\omega$ such that $K(U_\omega,Y)=\omega(Y)$ for any $Y \in \g$. Since $K$ is $ad(\g)$-invariant, $A(d\omega)$ coincides with the commutator of $U_\omega$, that is, $\{X \in \g, \ [X,U_\omega]=0\}.$ But the dimension of this subspace is greater than or equal to $r$. Then $n-c+1 \geq r$, that is, $c \leq n-r+1.$
\begin{corollary}
Let $\g$ be a semi-simple Lie algebra of rank $r$. There exists a form of maximal class $c=\dim \g$ if and only if $r=1$. \end{corollary}
\pf From the previous inequality, $c=n$ implies $r=1$. Thus $\g$ is of rank $1$ and it is isomorphic to $sl(2)$ in the complex case,  $sl(2,\R)$ or $so(3)$ in the real case. It is easy in each case to find a contact form. \begin{proposition}
Let $\g$ be a real $n$-dimensional Lie algebra such that for any $\omega  \in g^*$ with $\omega \neq 0$ we have $cl(\omega)=n$. Then $\g$ is of dimension $3$ and is isomorphic to $so(3)$.
\end{proposition}
\pf Let $I$ be a non trivial abelian ideal of $\g$. There exists $\omega \neq 0 \in \g^*$ such that $\omega (X)=0$ for any $X \in I$. This implies that $I \subset A(d\omega)$ and $cl(\omega) \leq n - dim I$ and this form is not of class $n$. Thus $\g$ is semi-simple. But a semi-simple Lie algebra with a form of maximal class is of rank $1$. It is isomorphic to $sl(2)$ or $so(3)$. In the first case, we have a form of class $2$. In the second case, we have a basis of $so(3)^*$ such that $d\omega_1=\omega_2\wedge \omega_3$, $d\omega_2=\omega_3\wedge \omega_1$, $d\omega_3=\omega_1\wedge \omega_2$. In this case any $\omega\neq 0 \in \g^*$ is of class $3$.
\medskip

\noindent{\bf Remark.} The  Proposition \ref{simple} give an upper bound to the set of the classes of linear form on a simple Lie algebra of rank $r$. In \cite{GozeCras1}, we found also a lower bound in case of $\g$ is simple of classical type. This bound is $2r$, 

\section{A contact form on $SL(2p)$}
Recall that $SL(2p)$ is a simple Lie group of rank $2p-1$ of dimension $4p^2-1$ whose elements are the matrices $M=(a_{i,j})$, $1 \leq i,j \leq 2p$ with determinant equal to $1$. It can be considered as an algebraic sffine variety  embedded in $\K^{4p^2}$ and defined by the polynomial equation $\det M=1$.  If we denote by $A_{i,j}$ the minor of $a_{i,j}$, then
$$\det M= \sum _{j=1}^{2p} (-1)^{i+j}a_{1,j}A_{1,j}.$$
We denote this determinant by $\Delta_M$.
\begin{lemma}
$$d\Delta_M= \sum (-1)^{i+j}A_{i,j}da_{i,j}.$$
\end{lemma}
Then a vector field $X$ on $\R^{4p^2}$ induces a vector field, denoted also $X$ on $SL(2p)$ if and only if $X(\Delta_M)=0$. We can construct a basis of left invariant vector fields. Let $L_M: M_1 \ra M\cdot M_1$ the left translation. We have
$$(L_M)^T(Id)= M \otimes Id$$
where $(L_M)^T(Id)$ is the tangent map to $L_M$ at the identity element $Id$ of $SL(2p)$. Likewise, if $R_M$ is the right translation by $M$, 
$$(R_M)^T(Id)= Id \otimes  ^tM.$$
If we denote by $\partial_{i,j}$ the derivate $\ds \frac{\delta}{\delta a_{i,j}}$, a basis of $sl(2p)$ is given by
$$\begin{array}{ll}
     X_{1,1} =& \ds\sum_{i=1}^{2p} a_{i,1}\partial_{i,1}  - \sum_{i=1}^{2p} a_{i,2p}\partial_{i,2p}\\
     \medskip
    X_{2,2} =&\ds \sum_{i=1}^{2p} a_{i,2}\partial_{i,2}  - \sum_{i=1}^{2p} a_{i,2p}\partial_{i,2p}\\
    \cdots & \cdots \\
    X_{2p-1,2p-1}=& \ds\sum_{i=1}^{2p} a_{i,2p-1}\partial_{i,2p-1}  - \sum_{i=1}^{2p} a_{i,2p}\partial_{i,2p}\\
    
\end{array}
$$
which form a basis of the Cartan subalgebra of $sl(2p)$ and by  the vectors
$$\begin{array}{ll}
     X_{1,2} =& \ds\sum_{i=1}^{2p} a_{i,1}\partial_{i,2} \\
     \medskip
    X_{2,1} =&\ds \sum_{i=1}^{2p} a_{i,2}\partial_{i,1}  \\
    \cdots & \cdots \\
    X_{k,l}=& \ds\sum_{i=1}^{2p} a_{i,k}\partial_{i,l}   \ \  k<l \\
     \medskip
 X_{l,k}=& \ds\sum_{i=1}^{2p} a_{i,l}\partial_{i,k}   \ \  k<l,\\

\end{array}
$$
We deduce a basis of $\g^*$:
\begin{equation}
\label{alpha}
\left\{
\begin{array}{ll}
\alpha_{k,k} = &\ds \sum_{i=1}^{2p} (-1)^{i+k}A_{i,k}d_{i,k}, \ \ i=1,\cdots 2p-1 \\
\alpha_{k,l}=& \ds\sum_{i=1}^{2p} (-1)^{i+k}A_{i,k}d_{i,l}   \ \  k<l\\
\alpha_{l,k}=& \ds\sum_{i=1}^{2p} (-1)^{i+l}A_{i,l}d_{i,k}   \ \  k<l\\
\end{array}
\right.
\end{equation}

The Cartan class of the linear form $\alpha_{1,1}+\alpha_{2,2}+\cdots +\alpha_{2p-1,2p-1}$ is $4p^2-2p+1$, that is the possible maximum. 

\noindent{\bf Remark.} The left invariance of these vector fields or of the linear forms is equivalent to $[X_{i,j},Y_{k,l}] =0$ or $\mathcal{L}(Y_{k,l})\alpha_{i,j}=0$ where $Y_{k,l}$ is a right invariant vector fields. A basis of these vector fields is given by
\begin{equation}\label{rvf}
\begin{array}{ll}
Y_{k,k}= &\ds\sum_{i=1}^{2p} a_{k,i}\partial_{i,1}  - \sum_{i=1}^{2p} a_{2p,i}\partial_{2p,i}\\
Y_{k,l}=& \ds\sum_{i=1}^{2p} a_{l,i}\partial_{k,i}   \ \  k<l\\
Y_{l,k}=& \ds\sum_{i=1}^{2p} a_{k,i}\partial_{l,i}   \ \  k<l\\
\end{array}
\end{equation}
\begin{theorem}
Let $\omega$ be the Pfaffian form
$$\omega= \sum_{i=1}^{2p}\sum_{j=1}^{p} a_{i,2j}da_{i,2j-1}-a_{i,2j-1}da_{i,2j}$$
Then $\omega$ induces a contact form on $SL(2p)$.
\end{theorem}
\pf We have 
$$d \omega= 2 \sum_{i=1}^{2p}\sum_{j=1}^{p} da_{i,2j}\wedge da_{i,2j-1}$$
If $\Theta$ means the "volume" form $(d \omega)^{2p^2}$ and $\Theta_{i,j}$ the $(4p^2-2)$-exterior form obtained from $\Theta$ from which the factor $ 
da_{i,2j}\wedge da_{i,2j-1}$ was removed, then
$$\omega \wedge (d\omega)^{2p^2-1}=2^{2p^2-1}\sum_{i=1}^{2p}\sum_{j=1}^{p} \Theta_{i,j} \wedge (a_{i,2j}da_{i,2j-1}-a_{i,2j-1}da_{i,2j}).$$
Then
$$
\begin{array}{ll}
\medskip
\omega \wedge (d\omega)^{2p^2-1}\wedge d\Delta_M=&2^{2p^2-1} \ds\sum_{i=1}^{2p}\sum_{j=1}^{p} ((-1)^{i+j}(\Theta_{i,j}\wedge a_{i,2j}da_{i,2j-1} \wedge A_{i,2j}da_{i,2j}\\ &
-(-1)^{i+j+1}\Theta_{i,j}\wedge (a_{i,2j-1}da_{i,2j} \wedge A_{i,2j-1}da_{i,2j-1})
\end{array}
$$
that is 
$$\omega \wedge (d\omega)^{2p^2-1}\wedge d\Delta_M=-2^{2p^2+p-1}\Delta_M\Theta$$
and $\omega$ induces a contact form on $SL(n)$.

\section{On the invariance of the form $\omega$}
\subsection{The Reeb vector field} Recall that  $\alpha$ is a contact form on a differential manifold $V$, then 
the
associated Reeb vector field is the section, $R_\alpha$, of $TM$ that generates the kernel of $d\alpha$ and
pairs with $\alpha$ to give $1$. Then it is defined by
$$\alpha(R_\alpha)=1, \ \ i(R_\alpha)dR_\alpha=0$$
where $i(R_\alpha)d\alpha (Y)= d\alpha (R_\alpha,Y)$. If we consider The Lie group $SL(n)$ and the contact form $\omega$, then 
\begin{proposition}
The Reeb vector field associated with the contact form $\omega$ on $SL(2p)$ is the vector field induced by
$$R_\omega=\ds \frac{1}{4\Delta} \sum_{i=1}^{2p}\sum_{j=1}^{p}(-1)^{i+1} ( A_{i,2j} \partial_{i,2j-1}+A_{i,2j-1}\partial_{i,2j})$$
\end{proposition}
\subsection{$J$-invariant forms}
Let $G$ be a Lie group and $J$ a linear subspace of the Lie algebra $\g$ of $G$.
\begin{definition}A $p$-form $\theta$ on $G$ is $J$-invariant if $L_{Y}\theta =0$ for any $Y\in J$.
\end{definition}
In particular, since $\g$ is the Lie algebra of left invariant vector fields on $G$, if $ J=\g$ i a $J$-invariant form is left invariant form on $G$. 

 In any point $g \in G$, $J$ defines a linear subspace $J_g$ of $T_gG$ and the distribution $\{J_g\}_{g \in G}$ is regular, that is $\dim J_g=\dim J$. If $\omega$ is a $J$-invariant Pfaffian form on $G$, the singular set $\Sigma_\omega$ associated with $\omega$ is the subset of $G$:
$$\Sigma_\omega=\{ g \in G, \ \omega_g(Y)_g)=0 \ \forall Y \in J\}.$$
\begin{lemma} (\cite{Lu1}). If $\omega$ is a $J$-invariant contact form on $G$, thus for any basis $\{Y_1,\cdots,Y_k\} $ of $J$, the map
$$\varphi(g)=(\omega(Y_1)(g),\cdots,\omega(Y_k)(g))$$
with values in $\R^k$ is of rank $k$ on $\Sigma$ and of rank $k-1$ on $G - \Sigma$ where $k=\dim J$.
\end{lemma}
In particular,we deduce
\begin{itemize}
  \item If $k \leq p+1$, then $\Sigma=\emptyset.$
  \item If $J$ is an abelian subalgebra of $\widetilde{\g}$, then $k \leq p+1.$
\end{itemize}
For example, there exists on the torus $\mathbb{T}_5$ a $J$-invariant contact form where $J$ is the abelian $2$-dimensional Lie algebra \cite{Lu1}.

\medskip

On a connected Lie group $G$, one defines classically the following complexes and their corresponding cohomologies:
\begin{enumerate}
  \item The de Rham cohomology of$G$, $H^*_{DR}(G)$, corresponding to the complex of differential forms on $G$.
  \item The cohomology of left invariant forms $H_L^*(G)$ corresponding to the complex of left invariant forms on $G$.
  \item The cohomology $H^*(G,M)$ of the group $G$ with values in a $G$-module $M$. The $p$-cochains are the mappings on $G^p$ with values in $M$.
  \item The differentiable cohomology $H_{diff}^*(G,M)$ where $M$ is a differentiable $G$-module. The $p$-cochains are the differentiable mappings on $G^p$ with values in $M$.
\end{enumerate}
and on the Lie algebra $\g$ of $G$
\begin{enumerate}
  \item The Chevalley-Eilenberg cohomology $H^*(\g,V)$ associated with the complex of exterior forms on $\g$ with values on a $\g$-module $V$.
  \item The $\h$-basic cohomology $H ^*(\g,\h,\K)$ where $\h$ is a Lie subalgebra of $\g$, and the $p$-cochains of the corresponding complex, the skew linear $p$-forms on $\g$ whose characteristic space contains  $\h$.
\end{enumerate}
  We have the following classical isomorphismes:
\begin{enumerate}
  \item $H_L^*(G)= H^*(\g,\R)$
  \item If $G$ is compact, $H^*_{DR}(G)=H_L^*(G)$
  \item If $G$ is semisimple and if $\h$ is the Lie algebra of the maximal compact subgroup of $G$, then $H^*_{diff}(G,\R)=H^*(\g,\h,\R)$.
\end{enumerate}
It is easy to define a complex on $G$ taking account the notion of $J$-invariance.
Let $\Lambda^p_J(G)$ the space of $J$-invariant $p$-form on the Lie group $G$. Since
$$L_X \circ d = d\circ L_X$$
for any vector field $X$, we have
$$\theta \in \Lambda^p_J(G) \Rightarrow d\theta \in \Lambda^{p+1}_J(G).$$
This permits to define a complex where the spaces of cochains are $\Lambda^p_J(G)$ and the coboundary operator the classical exterior differential. We denote by $(\Lambda^p_J(G),d)_p$ this complex:
$$
\Lambda^0_J(G) \stackrel{d} \ra \Lambda^1_J(G) \stackrel{d} \ra \Lambda^2_J(G) \ldots \stackrel{d}\ra \Lambda^n _J(G) \stackrel{d}\ra 0$$
where $n=\dim G$ and $\Lambda^0_J(G)$ is the space of functions on $G$ which are $J$-invariant. The corresponding cohomology $H^*_J(G)$ will be called the $J$-invariant de Rham cohomology. If $J=\widetilde{\g}$, thus
$H^*_J(G)$ coincides with the Chevalley-Eilenberg cohomology $H^*_{CE}(G)$ of $\g$ with values in the trivial module $R$. If $G$ is compcat, thus
$$H^*_J(G)=H^*_{dR}(G)=H^*_{CE}(\g,\R)$$
where $ H^*_{dR}(G)$ is the de Rham cohomology of $G$. These equalities are in general not satisfied when $G$ is not compact. For example, if $G=\R^2$ and $J$ the subalgebra generated by $\displaystyle \frac{\partial}{\partial x}$, thus $\Lambda^0_J(G)=\{f:R^2\rightarrow \R, \displaystyle \frac{\partial f}{\partial x}=0\}$ that is $\Lambda^0_J(G)=\{f(y)\}$, but $\Lambda^0_{dR}(G)=\{f(x,y)\}$, $\Lambda^0_{\R^2}(G)=\R.$

\medskip

\noindent{\bf Examples.}
\begin{enumerate} 
\item Contact form on the torus $\Bbb{T}^3$ invariant by $\Bbb{T}^1$
$$\omega= \cos(n_1\theta_1)d\theta_2 +\sin(n_1\theta_1)d\theta_3.$$ 
\item Contact form on the torus $\Bbb{T}^5$ invariant by $\Bbb{T}^2$ \cite{Lu1}.
$$
\begin{array}{ll}
\omega = & (\sin \theta_1  \cos \theta_3 - \sin \theta_2 \sin \theta_3) d\theta_4 + (\sin \theta_1 \sin \theta_3
+ \sin \theta_2 \cos\theta_3)d\theta_5 \\
&+ \sin \theta_2 \cos \theta_2 d\theta_1 - \sin\theta_1 \cos\theta_1 d\theta_2
+ \cos \theta_1 \cos \theta_2 d\theta_3
\end{array}
$$
which is invariant with respect to $\frac{\partial}{\partial \theta_4}$ and 
$\frac{\partial}{\partial \theta_5}.$
\item In \cite{G.R.JDGA} we give some contact forms on the $3$-dimensional Heisenberg Lie groups unvariant by a $2$-dimensional subgroup with an empty singular set $\Sigma.$
\end{enumerate} 
\subsection{Invariance of $\omega$}
Let us determine the symmetry group of the contact form $omega$ on $SL(2p)$. 
If $e$ denotes the unit of this Lie Group, then 
$$\omega (e)=\sum_i da_{i,i+1}(e)-da_{i+1,i}(e).$$
Let $a=(a_{i,j})$ be in $SL(2p)$. We have seen that $(L_a)^T(e)=a \otimes Id$, then 
$$(L_a^T)*^{-1}\omega (e)=\ds \sum_{i=1}^{2p}\sum_{j=1}^p (-1)^{i}(A_{i,2j}da_{i,2j-1}+A_{i,2j-1}da_{i,2j})$$
In fact $(L_a^T)*^{-1}= a^{-1} \otimes Id.$ This form coincides with $\omega (a)$ if and only if
$$
(-1)^iA_{i,2j}=a_{i,2j}, \ \ (-1)^iA_{i,2j-1}=-a_{i,2j-1}, \ \ i=1,\cdots 2p, \ \ j=1,\cdots,p.
$$
This shows that 
$$a=\begin{pmatrix}
    A_{1,1}  &  -A_{1,2} & A_{1,3} & \cdots & -A_{1,2p} \\
    \cdots &  \cdots &\cdots &\cdots &\cdots \\
    -A_{2p,1}  &  A_{2p,2} & -A_{2p,3} & \cdots & A_{2p,2p} 
\end{pmatrix}
$$
that is $$a=^ta^{-1}.$$
\begin{theorem}
The Pfaffian form $\omega=\ds \sum_{i=1}^{2p}\sum_{j=1}^{p} a_{i,2j}da_{i,2j-1}-a_{i,2j-1}da_{i,2j}$ induces on $SL(2p)$ a contact form which is invariant by the left action of the subgroup $SO(2p)$, and this invariance is maximal.
\end{theorem}
The maximality means that no other subgroup of $SL(2p)$ containing $SO(2p)$ left invariant the contact form $\omega$.

\pf If $a,b \in SL(2p)$, we have seen that the Jacobian matrix of the left translation at the identity is
$$(L_a)^T_e= a \otimes Id.$$
This implies
$$(^t(L_a)^T_e)^{-1}=(^ta)^{-1} \otimes Id.$$
Since 
$$\omega (e)= -da_{1,2}+da_{2,1}-d_{3,4}+da_{4,3} - \cdots -da_{2p-1,2p}+da_{2p,2p-1}$$
we obtain
$$(^t(L_a)^T_e)^{-1}\omega(e)=\sum_{i=1}{2p}(-1)^i\sum_{j=1}^{p}(A_{i,2j} da_{i,2j-1}+A_{i,2j-1}da_{i,2j}$$
where $A_{i,j}$ is the minor of $a$ associated with the coefficient $a_{i,j}$.
Assume that  $(^t(L_a)^T_e)^{-1}\omega(e)=\omega (a).$ Then
$$ a_{i,j}=(-1)^{i+j}A_{i,j}.$$
This is equivalent to 
$$^ta=a^{-1}$$
that is $a \in SO(2p).$

\medskip

\noindent{\bf Remark.} When $p=1$, that is in $SL(2)$, the form $\omega$ is induced by
$$\omega (a)=a_{1,2}da_{1,1}-a_{1,1}da_{1,2} +a_{2,2}da_{2,1}-a_{2,1}da_{2,2}.$$
Its decomposition in the basis (\ref{alpha}) is 
$$\omega (a) =u_{1,1}(a) \alpha_{1,1} + u_{1,2}(a) \alpha_{1,2} +u_{2,1}(a) \alpha_{2,1} $$
with
$$\begin{array}{l}
\medskip
    u_{1,1}(a)  =    \omega (a) (X_{1,1})= 2(a_{1,1}a_{1,2}+a_{2,1}a_{2,2})=2 <C_1,C_2>\\
    \medskip
    u_{1,2}(a)  =    \omega (a) (X_{1,2})= -a_{1,1}^2-a_{2,1}^2=- <C_1,C_1>\\ 
    \medskip
        u_{2,1}(a)  =    \omega (a) (X_{1,2})= a_{1,2}^2+a_{2,2}^2=<C_2,C_2>\\ 
\end{array}
$$
where $C_i$ is the column $i$ of the matrix $a$ and $<,>$ the classical inner product on the vector space of columns of square matrices. 

We can generalize this last decomposition for $SL(2p)$:
$$\omega (a)= \sum u(i,j)(a)\alpha_{i,j}.$$
To compute the coefficients $u(i,j)$ we remark that we can write $\omega$  of the following form:  If $d$ denotes the matrix $d=(d_{i,j})=da_{i,j}$ and $D_j$ the column $D_j=^t(d_{1,j}, \cdots,d_{2p,j})$ then $\omega(a)=\sum_{k=1}^{p} <C_{2k},D_{2k-1}>-<C_{2k-1},D_{2k}>$. Similary, the  left invariant vectors field $X_{l,k}$ can be written  
$$
\begin{array}{l}
 X_{k,k}=<C_k,\Delta_k>-<C_{2p},\Delta_{2p}>, \  \ \ k=1,\cdots,p  \\
    X_{k,l}=<C_k,\Delta_l>, \  \ \ 1\leq k < l \leq 2p \\
    X_{l,k}=<C_l,\Delta_k>,  \ \ \  1\leq k < l \leq 2p \\
\end{array}
$$ 
where $\Delta_k$ is the column $k$ of the matrix $(\partial_{i,j})$. Using these notations, we have
$$
\begin{array}{llll}
 u_{2k-1,2k-1}&=&<C_{2k-1},C_{2k}>+<C_{2p-1},C_{2p}>, & k=1,\cdots,p  \\
  u_{2k,2k}&=&-<C_{2k-1},C_{2k}>+<C_{2p-1},C_{2p}>, & k=1,\cdots,p  \\
    u_{k,2l}&=&-<C_k,C_{2l-1}>, & 1\leq k <2l \leq 2p \\
     u_{k,2l-1}&=&<C_k,C_{2l}>, & 1\leq k <2l -1\leq 2p \\
   u_{l,2k}&=&-<C_l,C_{2k-1}>,  & 1\leq 2k < l \leq 2p \\
     u_{l,2k-1}&=&<C_l,C_{2k}>,  &  1\leq 2k < l \leq 2p \\
\end{array}
$$ 

\subsection{The singular set $\Sigma_\omega$}
Recall that a Pfaffian form $\alpha$ on a Lie group $G$ is left invariant if and only if $\mathcal{L}(Y)(\alpha)=0$ for any right invariant vector field $Y$ on $G$, where $\mathcal{L}(Y)$ is the Lie derivate with respect the vector field $Y$ that is $\mathcal{L}(Y)(\alpha)=i(Y)d\alpha+d(\alpha(Y)).$ In the previous section, we have computed the right invariant vector fields on $SL(2p)$, they are linear combinations of the $Y_{k,l}$ defined in (\ref{rvf}). Recall also that, if $a \in SL(2p)$ and if $R_a$ is the right translation by $a$ in this group, $R_a(b)=ba$ then its tangent map at the identity is
$$(R_a)T_e=Id \otimes ^t a.$$
Its dual map is
$$(^t(R_a)^T_e)^{-1}= Id \otimes a^{-1}.$$
Since $\omega(e)=\sum_{i=1}^p (-da_{2i-1,2i}+da_{2i,2i-1})$ then $(^t(R_a)^T_e)^{-1}\omega (e)=\omega (a)$ if and only if
 $$\omega (a)=\sum_{i=1}^p\sum_{j=1}^{2p}(-1)^{j-1} (A_{2i,j}da_{2i-1,j}+A_{2i-1,j}da_{2i,j})$$
 giving
$$a=J_{2p}^{-1}a^{-1} J_{2p}$$
where$$J_{2p}=Id_{p} \otimes J_2$$ with $$J_2=\begin{pmatrix}
   0   & 1  \\
  -1    &  0
\end{pmatrix}
$$ 
\begin{theorem}
The Pfaffian form $\omega=\ds \sum_{i=1}^{2p}\sum_{j=1}^{p} a_{i,2j}da_{i,2j-1}-a_{i,2j-1}da_{i,2j}$ induces on $SL(2p)$ a contact form which is invariant by the right action of the subgroup $H$ where
$$H=\{a \in Sl(2p), \ \ / \ ^t aJ_{2p}a=J_{2p}\}$$ and this invariance is maximal.
\end{theorem}
It is clear that $H$ is a closed Lie subgroup of $SL(2p)$. In fact, if $a,b \in H$,
$$^t(ab)J(ab)=^tb^taJab=^tbJb=J$$
and $ab \in H$. Likewise
$$^t(a)^{-1}Ja^{-1}=(-aJ^ta)^{-1}=^tJ^{-1}=J.$$
The Lie algebra $\h$ of $H$ is the set of square matrices $Y$ which satisfy
$$JY+^tYJ=0.$$

For example, if $p=1$ the elements of $\h$ are the matrices which satisfy
$$\begin{pmatrix}
    a_{1,1} & a_{2,1}   \\
    a_{1,2}& a_{2,2}\end{pmatrix} \cdot \begin{pmatrix}
    0  & -1   \\
    1  & 0 
\end{pmatrix} + 
\begin{pmatrix}
    0  & -1   \\
    1  & 0 
\end{pmatrix}\cdot\begin{pmatrix}
    a_{1,1} & a_{1,2}   \\
    a_{2,1}& a_{2,2}\end{pmatrix} = \begin{pmatrix}
  0    &  - a_{1,1} - a_{2,2}   \\
    a_{1,1} + a_{2,2}   &  0
\end{pmatrix} =0
$$
and $H=SL(2)$. This shows that $\omega$ induces on $SL(2)$ a right- invariant contact form. Then the singular set $\Sigma_\omega =\emptyset.$

If $p=2$, then $\h$ is the Lie subalgebra of $sl(3)$ defined by
$$\h=\left\{ \begin{pmatrix}
   a_{1,1} & a_{1,2} &a_{1,3} & a_{1,4} \\
   a_{2,1} & -a_{1,1} &a_{2,3} & a_{2,4} \\
  - a_{2,4} & a_{1,4} &a_{3,3} & a_{3,4} \\
    a_{2,3} & -a_{1,3} &a_{4,3} & -a_{3,3} \\
\end{pmatrix}  \right\}
$$
and $H$ is a $10$-dimensional Lie group isomorphic to $Sp(2)$. In the general case, $\h$ is the set of matrices
$(M_{i,j})$ avec $1 \leq i,j \leq p$ avec
$$M_{i,i}=\begin{pmatrix}
  a_{2i-1,2i-1}    &  a_{2i-1,2i}    \\
    a_{2i,2i-1}    &   -a_{2i-1,2i-1} 
\end{pmatrix}, \ \ M_{i,j}=\begin{pmatrix}
     a_{2i-1,2j-1}   &  - a_{2i-1,2j}  \\
     a_{2i,2j-1}   &   a_{2i,2i} 
\end{pmatrix},  \ \  M_{j,i}=JM_{i,j}J \ \ {\rm for} \ i<j.$$
Its dimension is $p(2p+1).$

\subsection{The Lie group $SO(n)$}

Recall that $SO(n)$ is the simple Lie group whose elements are the orthogonal matrices of determinant $1$.  Its Lie algebra, $so(n)$, is the vector space of skew-symmetric matrices. We have seen, in the previous section, that the contact form $\omega$ on $SL(2p)$ is left invariant by the action of $SO(2p)$. which is of dimension $p(2p-1)$. Let $a=(a_{i,j})$ be in $S0(n)$. This is equivalent to say that the coefficients $a_{i,j}$ are solutions of the algebraic system
$$ f_{k,l}=\sum_j a_{j,k}a_{j,l}-\delta_k^l=0, \ \ 1 \leq k \leq l \leq n.$$
In particular a differential form $\phi$ defined on $GL(n,\R)$ induces a non trivial form on $SO(n)$ if and only if
$$\phi \wedge df_{1,1} \wedge df_{1,2} \wedge \cdots \wedge df_{n,n} \neq 0$$
on $SO(n)$. 

As before, we denote by $d_{i,j}$ the differential $da_{i,j}$. We consider on the set of couples $(i,j)$ the lexicographic order. Then
$$df_{1,1} \wedge df_{1,2} \wedge \cdots \wedge df_{n,n} =\sum M_{I} d_{i_1,j_1} \wedge d_{i_2,j_2} \wedge \cdots \wedge d_{i_N,j_N}$$
with $(1,1) \leq (i_1,j_1) < \cdots < (i_N,j_N) \leq (n,n), \ N=\ds \frac{n(n+1)}{2}$and $I=\{(i_1,j_1),\cdots,(i_N,j_N)\}.$ To compute the coefficients $M_I$, we consider the matrix $M$ of $n^2$ lines and $N$ columns which represents the coefficients of the ordered free family $\{1/2df_{1,1}, f_{1,2},\cdots,df_{n-1,n},1/2f_{n,n}\}$ in the ordered basis $d_{i,j}$. Then $M_I$ is the minor of order $N$ whose lines are indexed by the coefficients $(i,j) \in I$. 

If $n=3$, $N=6$, we consider the map $F=(1/2f_{1,1},f_{1,2}, \cdots, 1/2f_{n,n}): GL(3) \ra \R^6$ whose components are written in the lexicographic order. Its jacobian matrix is of order $6 \times 9$. Then the coefficient $M_I$ is the determinant of the restricted matrix of order $6$ correspondant to the sequence $I$. Since the matrix $(a_{i,j} )\in SO(3)$, each one of these coefficients are $0$ or of type $\pm a_{i,j}a_{k,l}.$

\medskip

Let us return to the general case. We have seen that the tangent map to the left translation by $a \in GL(n)$ is
$$(L_a) ^T_E=a \otimes Id.$$
If $a \in SO(n)$, then $a^{-1}=^ta$ and its dual map ${(L_a)^T_E}^*=a \otimes Id$. To simplify, we will write this map $l_a$. A Pfaffian form $\omega$ on $SO(n)$ is left invariant if $l_a(\omega(e))=\omega(a)$. We consider on $GL(n)$ the following global basis of Pfaffian forms, written in this order,
$$\{d_{11},d_{12},\cdots,d_{1n},d_{21},\cdots,d_{2n},\cdots,d_{n,1},\cdots,d_{nn}\}$$
and the vector $e_{i,j}=(0,\cdots,1,0,\cdots,0)=d_{ij}$.  With these notations we have
$$df_{1,1}(a)=l_a(e_{1,1}),df_{2,2}(a)=l_a(e_{2,2}),\cdots,df_{n,n}(a)=l_a(e_{n,n})$$
and for $i<j$,
$$df_{i,j}(a)=l_a(e_{i,j}+e_{j,i}). $$
Once again to simplify notations, we put $l_{i,j}=l_a(e_{i,j})$. If $\Theta_n$ is the $\frac{n(n+1)}{2}$-exterior form
$$\Theta_n=\wedge_{1 \leq i \leq j \leq n} df_{i,j}$$
then
$$\Theta_n=\wedge_{1 \leq i \leq n}l_{i,i}\wedge_{1 \leq i< \leq n} (l_{i,j}+l_{j,i}).$$

\bigskip

For example, if $n=3$ and $a=(a_{ij})$ be an element of $SO(3).$ The Pfaffian form $\omega$ (or rather the Pfaffian form induced by $\omega$ on $SO(3)$) given by $\omega=a_{11}d_{12}+a_{21}d_{12}+a_{31}d_{32}$ is left-invariant.The coefficients of $\Theta_3$ are the minors of order $6$ of the matrix
$$\begin{pmatrix}
     a_{11} & 0 & 0 & a_{21} & 0 & 0 & a_{31} & 0 & 0\\
     a_{12} &  a_{11} & 0 & a_{22} & a_{21} & 0 & a_{32} & a_{31} & 0\\
     a_{13} & 0 & a_{11} & a_{23} & 0 & a_{21} & a_{33} & 0 & a_{31} \\
     0 & a_{12} & 0 & 0 & a_{22} & 0 & 0 & a_{32} & 0\\
     0 & a_{13}& a_{12} & 0 & a_{23} & a_{22} & 0 & a_{33} & a_{32} \\
     0 & 0 & a_{13} & 0 & 0 & a_{23} & 0 & 0 &a_{33} 
\end{pmatrix}
$$
If we denote by $M_{i,j,k}$ this minor obtained after removing the colums $i,j,k$, then 
$$\omega \we d\omega \we \Theta_3= 
a_{11}(M_{2,4,5}+M_{2,7,8})-a_{21}(M_{1,2,5}+M_{5,7,8})+a_{31}(M_{1,2,8}+M_{4,5,8})V$$
where $V$ is the volume form $V=d_{11} \we d_{12} \we \cdots \we d_{32} \we d_{33}$.
But
$$
\begin{array}{ll}
M_{2,4,5}=a_{32}(a_{12} a_{31} - a_{11} a_{32}) \Delta=a_{32}a_{23} & M_{2,7,8}=-a_{22}(a_{11}a_{22}-a_{12}a_{21})\Delta=-a_{22}a_{33}\\
M_{1,2,5}=a_{32}(a_{21}a_{32}-a_{22}a_{31})\Delta=a_{32}a_{13} & M_{5,7,8}=-a_{12}(a_{11}a_{22}-a_{12}a_{21})\Delta=-a_{12}a_{33}\\
M_{1,2,8}=a_{22}(a_{21}a_{32}-a_{22}a_{31})\Delta=a_{22}a_{13} & M_{4,5,8}=a_{12}(a_{11}a_{32}-a_{12}a_{31})\delta=-a_{12}a_{23}
\end{array}
$$
So
$$\omega \we d\omega \we \Theta_3= a_{11}(a_{32}a_{23}-a_{22}a_{33})-a_{21}(a_{32}a_{13}-a_{12}a_{33})+a_{31}(a_{22}a_{13}-a_{12}a_{23})V=-\Delta V$$
and $\omega$ induces a contact form on $SO(3)$. In fact we find the property specifying that all non-zero left-invariant linear forms on this group are  contact forms. The general case will be studied later.

\end{document}